\documentclass [11pt]{article}
\usepackage{amsthm}
\usepackage{amsmath}
\usepackage{amssymb}
\usepackage{epsfig}

\newcounter{contador}
\newtheorem{propo}[contador]{Proposition}
\newtheorem{teo}[contador]{Theorem}
\newtheorem{lem}[contador]{Lemma}

\newtheorem{corol}[contador]{Corollary}

\newcommand{\rec}{\noindent}    
\newcommand{\dem}{\rec {\bf Proof. }}  
\renewcommand{\qed}{\ \hfill\rule[-1mm]{2mm}{3.2mm}\newline} 

\newcommand{\dps}{\displaystyle} 

\newcommand{\G}{\Gamma}

\newcommand{\su}{{\mathbb S}^1}  
\newcommand{\enya}{${\rm \tilde{n}}$}
\newcommand{\bx}{{\mathbf x}}

\newcommand{\sg}{{\cal G}}

\newcommand{\R}{{\mathbb R}}

\newcommand{\U}{{\cal{U}}}

\title{Some properties of the  $k$-dimensional Lyness' map}

\author{Anna Cima$^{(1)}$, Armengol Gasull$^{(1)}$ and V\'{\i}ctor Ma\~{n}osa $^{(2)}$
  \\*[.1truecm]
{\small \textsl{$^{(1)}$ Dept. de Matem\`{a}tiques, Facultat de
Ci\`{e}ncies,}}
\\*[-.25truecm] {\small \textsl{Universitat Aut\`{o}noma de Barcelona,}}
\\*[-.25truecm] {\small \textsl{08193 Bellaterra, Barcelona, Spain}}
\\*[-.25truecm] {\small \textsl{cima@mat.uab.cat, gasull@mat.uab.cat}}
\\*[-.25truecm]
\\*[-.25truecm] {\small \textsl{$^{(2)}$ Dept. de Matem\`{a}tica Aplicada III (MA3),}}
\\*[-.25truecm] {\small \textsl{Control, Dynamics and Applications Group (CoDALab)}}
\\*[-.25truecm] {\small \textsl{Universitat Polit\`{e}cnica de Catalunya (UPC)}}
\\*[-.25truecm] {\small \textsl{Colom 1, 08222 Terrassa, Spain}}
\\*[-.25truecm] {\small \textsl{victor.manosa@upc.edu}}}

\begin{document}
\maketitle

\begin{abstract}
 This paper is devoted to study some properties of the
 $k$-dimensional Lyness' map
$ F(x_1,\ldots,x_k)=(x_2,\ldots,x_k,(a+\sum_{i=2}^{k}x_i)/{x_1}).$
Our main result presentes a rational vector field that gives  a Lie
symmetry for $F.$ This vector field is used, for $k\le 5,$ to give
information about the nature of the invariant sets under $F.$ When
$k$ is odd, we also present a new (as far as we know) first integral
for $F\!\circ F$ which allows to deduce in a very simple way several
properties of the dynamical system generated by $F.$  In particular
for this case we prove that, except on a given codimension one
algebraic set, none of the positive initial conditions can  be a
periodic point of odd period.
\end{abstract}

\rec {\sl 2000 Mathematics Subject Classification} 39A20, 37E35.

\rec {\sl PACS numbers} 02.30.Ik, 05.45.-a.

 \rec {\sl Keywords:} Lyness' difference equations, Lie symmetry, first integral,
 integrable and non--integrable discrete dynamical system.\newline

\section{Introduction and main results}
The second and more specially the third order Lyness' difference
equations
$$
y_{n+2}=\frac{a+y_{n+1}}{y_{n}}\quad \mbox{ and }\quad
y_{n+3}=\frac{a+y_{n+1}+y_{n+2}}{y_{n}},\quad \mbox{ with}\quad a\geq 0,
$$
have been considered as emblematic examples of integrable discrete
systems, see for instance \cite{HKY,K,MT,RQ}. The dynamics of the
above equations, or their associated maps, has been the objective of
recent intensive investigation. Nowadays the behaviour of  the
orbits  is well known when positive initial conditions are
considered, see \cite{BR1,BC,CGM07-1,BMG,Z}.  However few results
have been obtained for negative initial conditions, see
\cite{CGM1,BMG}.

On the other hand, the study of  the  $2$ and $3$-dimensional
Lyness' maps is the starting point for the study of other birational
integrable maps, see for instance
\cite{BR2,CGM07-2,BMG,GM,I1,I2,IR1,RQ}, and \cite{CM} for a general
paper on this topic.

For $k\geq 4,$ very few results, apart of the ones obtained recently
by Bastien and Rogalsky \cite{BR}, are known for the $k$-th order
Lyness' equation

\begin{equation}\label{lynesseed}
y_{n+k}=\frac{a+\sum_{i=1}^{k-1}y_{n+i}}{y_n}.
\end{equation}

The main difference between the  $k=2,3$ and the $k\geq 4$ scenarios
is that the first cases are \textsl{integrable} in the sense that
the associated maps have $1$ and $2$ functionally independent first
integrals, respectively (in this paper we say that a map $F$ is
integrable if it has   $k-1$ functionally independent first
integrals). It seems that this property is not shared for the
Lyness' equations when $k\ge 4.$

Indeed, consider the  $k$-dimensional Lyness' map associated to the
difference equation~(\ref{lynesseed}),
\begin{equation}\label{lynk}
F(x_1,\ldots,x_k)=\dps{\left(x_2,\ldots,x_k,\frac{a+\sum_{i=2}^{k}x_i}{x_1}\right)} ,\mbox{ with}\quad a\geq 0,
\end{equation}
which is a diffeomorphism from ${\mathcal Q}^+:=\{{\bf
x}=(x_1,\ldots,x_k)\in\R^k \,:\, x_1>0, x_2>0, \ldots, x_k>0 \}$
into itself. It is well known that it has the following couple of
functionally independent first integrals
\begin{equation}\label{v1}
V_1({\bf x})=\left(a+\sum_{i=1}^k
x_i\right)\left(\prod_{i=1}^{k}(x_i+1)\right)/(x_1\cdots x_k)
\end{equation}
and
\begin{equation}\label{v2}
V_2({\bf x})=\left(a+\sum_{i=1}^k
x_i+x_1x_k\right)\left(\prod_{i=1}^{k-1}(1+x_i+x_{i+1})\right)/(x_1\cdots
x_k).
\end{equation}
A third functionally independent first integral for $k\geq 5$ has
recently been given in \cite{G}. Moreover in that paper it is
conjectured that for any $k,$ the Lyness' map has up to
$E\left(\frac{k+1}{2}\right)$ functionally independent first
integrals, where $E(\cdot)$ denotes the integer part function. The
conjecture seems to be true for $k\le6,$ see again \cite{G}.

The integrable structure for $k=2,3$ implies that  the dynamics of
the maps studied in the above references is in fact one-dimensional.
In any case, although  for $k\ge 4$ the above assertion seems not to
be true, the existence of several first integrals  reduces the
dimension of the space where the dynamics takes place.  If the above
conjecture would be true, then $k-E\left(\frac{k+1}{2}\right),$
would be, generically, the dimension of the invariant manifold given
by the level sets of the first integrals.

One  geometrical object that has played a key role to understand the
dynamics of a large class of 2 and $3$-dimensional maps is the
\textsl{Lie symmetry} of the map, \cite{CGM07-1,CGM07-2}. Recall,
that a vector field $\bf X$ is said to be a Lie symmetry of a map
$G$ if it satisfies the condition
\begin{equation}\label{estrella}
{\bf X}(G({\bf x}))=(DG({\bf x}))\,{\bf X}({\bf x}).
\end{equation}
The vector field $\bf X$ is related with the dynamics of $G$ in the
following sense: $G$ maps any orbit of the differential system
determined by the vector field, to another orbit of this  system,
see \cite{CGM07-2}. In the integrable case, where the dynamics are
in fact one dimensional, the existence of a Lie symmetry  fully
characterizes the dynamics. In \cite[Thm. 1]{CGM07-2} it is proved
that if $G:\U\to \U$ is a diffeomorphism  having a Lie symmetry $\bf
X$, and such $G$ preserves $\gamma,$ a solution of the differential
equation $\dot x={\bf X}(x),$ then the dynamics of
$\left.G\right|_{\gamma}$ is either conjugated to a rotation,
conjugated to a translation of the line, or constant, according
whether $\gamma $ is homeomorphic to $\su$, $\R$, or a point,
respectively. Other properties of the Lie symmetries of discrete
systems are studied in \cite{HBQC}.

One of the main results of this paper, is the following theorem
where the expression of a Lie symmetry for the $k$-dimensional
Lyness' map is given.

\begin{teo}\label{propoSL}
For $k\ge3,$ the  vector field ${\bf X}_k=\sum_{i=i}^{k}
X_i\dps{\frac{\partial }{\partial x_i}},$  is a Lie symmetry for the
$k$-dimensional Lyness' map (\ref{lynk}), where
\begin{equation}\label{X1}
X_1({\bf x})=\dps{\frac{(x_1+1)\left[\prod_{i=2}^{k-1}
(1+x_i+x_{i+1})\right](a+\sum_{i=1}^{k-1} x_i-x_2x_k)}{\prod_{i=2}^k
x_i}},\end{equation}
\begin{equation}\label{Xm}
X_m({\bf x})=\dps{\frac{(x_m+1)\left[\prod_{i=1,i\neq m-1,m}^{k-1}
(1+x_i+x_{i+1})\right](a+\sum_{i=1}^{k}
x_i+x_1x_k)(x_{m-1}-x_{m+1})}{\prod_{i=1,i\neq m}^k x_i}},
\end{equation}
for all  $2\leq m\leq k-1$, and
\begin{equation}\label{Xk}
X_k({\bf x})=-\dps{\frac{(x_k+1)\left[\prod_{i=1}^{k-2}
(1+x_i+x_{i+1})\right](a+\sum_{i=2}^{k}
x_i-x_1x_{k-1})}{\prod_{i=1}^{k-1} x_i}}.
\end{equation}
\end{teo}

Once we have the candidate ${\bf X}_k$ to be a Lie symmetry of the
Lyness' map $F$ the proof of Theorem \ref{propoSL} only will consist
in checking that (\ref{estrella}) holds. We give here some hints of
how we have found the above ${\bf X}_k.$ Observe that if there
exists a vector field $ {\bf X}_{k},$ satisfying equation
(\ref{estrella}) for the $k$-dimensional Lyness' map (\ref{lynk}),
then
$$
\left(\begin{array}{c}
  X_1(F) \\
  X_2(F) \\
  \vdots \\
  X_k(F)
\end{array}\right)=
\left(\begin{array}{cccccc}
  0 & 1 & 0 & 0 &\cdots & 0 \\
  0 & 0 & 1 & 0 & \cdots & 0 \\
  \vdots & {} & {} & {} & {} & {} \\
  -\frac{a+\sum_{i=2}^{k}x_i}{x_1^2} & \frac{1}{x_1} & {} & \cdots& & \frac{1}{x_1}
\end{array}\right)\left(\begin{array}{c}
  X_1 \\
  X_2 \\
  \vdots \\
  X_k
\end{array}\right).
$$
Hence it is necessary that
\begin{equation}\label{recuobt2}
X_{i+1} = X_i(F), \mbox{ for } i=1,\ldots,k-1,
\end{equation}
and the \emph{``compatibility condition''}:
\begin{equation*}\label{comp2}
 X_k(F)
=-\dps{\left(\frac{a+\sum_{i=2}^{k}x_i}{x_1^2}\right)}X_1+\dps{\frac{1}{x_1}\left[\sum_{i=2}^k
X_i \right]}.
\end{equation*}
Thus, the construction of a Lie symmetry is straightforward once the
right expression of $X_1,$  as a seed of equations (\ref{recuobt2}),
is obtained. In \cite{BC,CGM07-1} and \cite{CGM07-2} the expressions
of ${\bf X}_2$ and ${\bf X}_3$ are given. The idea of these papers
is that these vector fields have to be multiples of $\nabla V_1$ and
$\nabla V_1\times \nabla V_2,$ respectively. These constructions are
used to force  $F$ and ${\bf X}_k$ to share the same set of first
integrals. When $k\ge4$ we can not use anymore this idea because
there are no enough first integrals. Nevertheless we observe that
the first components of ${\bf X}_2$ and ${\bf X}_3,$ are
$$
\frac{(x_1+1)(a+x_1-x_2^2)}{x_2}\quad\mbox {and}\quad
\frac{(x_1+1)(1+x_2+x_3)(a+x_1+x_2-x_2x_3)}{x_2x_3},
$$
respectively. Thus it seems natural to try with $X_1$ as the
expression given in (\ref{X1}) and, indeed, it works!. From this
starting point, the proof for a given small $k$ is only a matter of
computations, while the proof for a general $k$ is long and tedious,
but straightforward. It is done in Section \ref{provaproposl}. We
suggest to skip this section in a first reading of the paper.

Another result that helps for understanding the dynamics generated
by $F,$ when $k$ is odd, is given in next proposition. In this
result, the key point is the existence of a new (as far as we know)
first integral for $F^2=F\!\circ F$ for any odd $k\ge3.$ As we will
see, our proof of the existence of this function is inspired in the
paper \cite{MT}, where this first integral is given for $k=3.$

\begin{teo}\label{propnova}Set $k=2\ell+1,$ ${\bf x}=(x_1,\ldots,x_{2\ell+1})$
and consider $F$ from $\mathcal{Q}^+$ into itself.
 Then

\begin{enumerate}

\item[(a)] The function

\begin{equation}\label {W} W({\bf x})=\frac{\prod_{j=0}^
{\ell}(x_{2j+1}+1)}{\prod_{j=1}^{\ell}x_{2j}},
\end{equation}
is a first integral of $F^2.$

\item[(b)] For any $\ell\geq 2$, the function $V_3:=W+W(F),$ which
is,
\begin{equation}\label{v3}
V_3({\bf x})=\frac{\prod_{j=0}^
{\ell}x_{2j+1}(x_{2j+1}+1)+(a+\sum_{j=1}^{2\ell +1}
x_j)\prod_{j=1}^{\ell}x_{2j}(x_{2j}+1)}{\prod_{i=1}^{2\ell+1}
x_{i}},
\end{equation}
is a first integral of $F$ which is functionally independent with
the first integrals   $V_1$ and $V_2$ given in (\ref{v1}) and
(\ref{v2}), respectively.

\item[(c)] The function $W\cdot W(F)$ coincides
with the first integral of $F,$ $V_1$ given in (\ref{v1}). In other
words,  $V_1=W\cdot W(F).$

\item [(d)] The algebraic set $\mathcal{G}:=\{{\bf x}\in{\cal Q}^+ \,:\,
W({\bf x})-W(F({\bf x}))=0\} $ is invariant by $F.$

\item[(e)] If the map $F$ has some periodic point of odd period then
it has to be contained in $\mathcal{G}.$

\item[(f)] Setting $\mathcal{G}^{\pm}:=\{{\bf x}\in{\cal Q}^+ \,:\, \pm(W({\bf
x})-W(F({\bf x})))>0\},$ the map $F$ sends ${\mathcal
 G}^{+}$ into ${\mathcal G}^{-}$ and viceversa, and both sets
 are invariant by $F^2.$ Furthermore, the dynamics of $F^2$ on each of these sets are conjugated,
being the map $F$ itself the conjugation.

\item[(g)] The measure
$$m_1(B):=\int_{B} \frac{\pm1}{\Pi({\bf x}) (W({\bf x})-W(F({\bf x}))) }\,d{\bf x},$$
where $\Pi({\bf x})=\prod_{i=1}^k x_i,$ is an invariant measure for
$F^2,$ {\it i.e.} $m_1(F^2(B))=m_1(B),$ where $B$ is any measurable
set in $\mathcal{G}^{\pm}.$

\item[(h)] The measure
$$m_2(B):=\int_{B} \frac{1}{\Pi({\bf x})}\,d{\bf x}, $$
is an  invariant measure for $F^2,$ {\it i.e.} $m_2(F^2(B))=m_2(B),$
where $B$ is any measurable set in ${\mathcal Q}^+.$

\end{enumerate}

\end{teo}

We remark that,  when $k$ is odd, the first integral $V_3$ given
above coincides with the one given recently in \cite{G}. Observe
also that the invariant algebraic surface $\mathcal G$ was already
found in \cite{CGM07-1}, but only  for $k=3.$ Also, as we will see
in Subsection \ref{dinreduida}, the function $W$ is  useful to make
an explicit simple order reduction when we study the dynamics of $F$
for $k$ odd.

Although, by using both theorems  we have not been able to present a
complete study of the higher dimensional Lyness' map, in next
results we give some information about the invariant sets in the
phase space when $k=4,5.$ We prove:

\begin{propo}\label{propoL4}
The vector field ${\bf X}_{4}$ given by  equations
(\ref{X1})--(\ref{Xk}) for $k=4$, is a Lie symmetry for the
$4$-dimensional Lyness' map. Moreover, ${\bf X}_4(V_i)=0$, for
$i=1,2,$ and then  the sets $I_{h,k}:=\{V_1=h\}\cap\{V_2=k\}\cap
{\cal Q}^+$ are invariant by $F$ and by the flow of ${\bf X}_{4}.$

Furthermore, if we assume that both first integrals intersect
transversally on $C_{h,k}$, a connected component of $I_{h,k},$ then
$C_{h,k}$ is diffeomorphic to a torus.
\end{propo}

\begin{propo}\label{propoL5} The vector field ${\bf X}_{5}$ given by
equations (\ref{X1})--(\ref{Xk}) for $k=5$ is a Lie symmetry for the
$5$-dimensional Lyness' map. Moreover,  ${\bf X}_5(V_i)=0$, for
  $i=1,2,3$, where
      $$\begin{array}{rl}
V_3({\bf x})= \displaystyle{\frac1{x_1 x_2 x_3 x_4 x_5}}\Big(&
x_1\, x_3\,
x_5\, \left( 1+x_1 \right) \left( 1+x_3 \right)  \left( 1+x_5\right)
\\&+   x_2\,x_4\,\left( 1+x_2
\right) \left( 1+x_4 \right)\left(
a+x_1+x_2+x_3+x_4+x_5 \right)\Big),
    \end{array}
  $$
and  the sets $I_{h,k,l}:=\{V_1=h\}\cap\{V_2=k\}\cap\{V_3=\ell\}\cap
{\cal Q}^+,$ which generically have  at least two connected
components, are invariant by $F$ and by the flow of ${\bf X}_{5}.$

Furthermore, if we assume that the three first integrals intersect
transversally on $C_{h,k,\ell}$, a connected component of
$I_{h,k,\ell},$ then $C_{h,k,\ell}$ is diffeomorphic to a
(two--dimensional) torus.
\end{propo}

It is important to notice that in the above results we do not assert
that most of the connected components   $C_{h,k}$  and
$C_{h,k,\ell}$ are tori, since we did not succeed to prove that over
them the intersection of the energy levels at the whole sets
$C_{h,k}$ and $C_{h,k,\ell}$  are transversal. Hence it remains open
to decide whether they are two--dimensional differentiable manifolds
or not, and in the case that they are not differentiable manifolds
to decide which are their topology. In any case, our result reduces
the problem to a computational question.

\centerline{\includegraphics[scale=0.55]{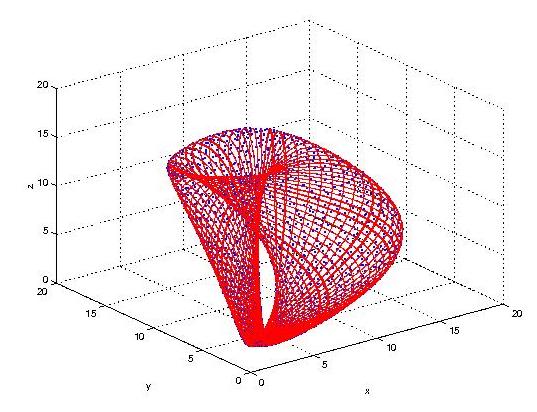}}
\begin{center}Figure 1. Projections into  $\R^3$ of the flow associated to
the Lie symmetry ${\bf X}_4$, and the orbit of the Lyness' map, for
$k=4$ and $a=4,$ both with initial condition (1,2,3,4).
\end{center}

Our numerical simulations seem to indicate that for $k=4,$ all the
generic level curves $I_{h,k},$ are connected. On the other hand,
for $k=5,$ generically the sets $I_{h,k,l}$ seem to have  exactly
two connected components. In Figures 1 and 2 we give a projection in
$\mathbb{R}^3$ of these surfaces. Indeed, in Figure 1 we represent,
both an orbit of $F$ and an orbit of ${\bf X}_4$ starting with the
same initial condition, and in Figure 2 an orbit of $F$ for $k=5.$
Notice that in both cases the behavior of the orbits seems to
indicate that two, (respectively three) is the maximum number of
independent first integrals for $F$ when $k=4$ (resp. $k=5$), as it
is suggested in \cite{G}.

\centerline{\includegraphics[scale=0.55]{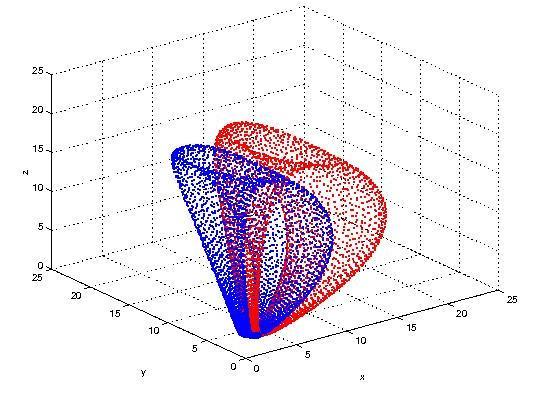}}
\begin{center}
Figure 2. Projection into  $\R^3$ of the first $5000$ iterates of
the Lyness' map for $k=5$ and  $a=1,$ staring at (1,2,3,4,5). Odd
and even iterates are in different connected components.\end{center}

A remarkable fact that  Figure 1 shows is that, although for $k=4$
the manifold $I_{h,k}$ is invariant for both the map $F$ and the
flow of ${\bf X}_4,$ the map $F$ seems to send an orbit of ${\bf
X}_4$ to a different orbit of the vector field. Further numerical
experiments seems to confirm this fact. Under this situation, we can
not use the techniques developed in \cite{CGM07-2}. This fact makes
more difficult the knowledge of the behavior of $F$ restricted to
each $I_{h,k}$ and is one of the important differences between the
cases $k=2,3$ and $k=4.$

\section{ Proof of Theorem \ref{propoSL}}\label{provaproposl}

As we have already explained in the introduction, the existence of a
vector field ${\bf X}_{k}=\sum_{i=i}^{k} X_i{\frac{\partial
}{\partial x_i}}$ satisfying equation (\ref{estrella}) for the
$k$-dimensional Lyness' map (\ref{lynk}),  is equivalent to the set
of equations
\begin{equation}\label{recuobt}
X_{i+1} = X_i(F), \mbox{ for } i=1,\ldots,k-1,
\end{equation}
together with the compatibility condition:
\begin{equation}\label{comp}
 X_k(F)
=-\dps{\left(\frac{a+\sum_{i=2}^{k}x_i}{x_1^2}\right)}X_1+\dps{\frac{1}{x_1}\left[\sum_{i=2}^k
X_i \right]}.
\end{equation}

The proof will consist in checking that the choice of ${\bf X}_k$
given in the statement satisfies equations (\ref{recuobt}) and
(\ref{comp}). The result is straightforward for $k=3,4,5$ and we
omit the details. So, from now on, we assume that $k\ge6.$

We proceed  in two steps:

\textsl{First step:} We will show that from the expression
(\ref{X1}) of $X_1$ as a seed of equations (\ref{recuobt}) we obtain
the expressions of $X_m$ for $m=2,\ldots,$$k-1$ and $X_k$ given by
equations (\ref{Xm}) and (\ref{Xk}), respectively.

\textsl{Second step:} We will prove that the compatibility condition
(\ref{comp}) is satisfied.

\vspace{0.4cm}

\noindent \textbf{First step:} We start with some preliminary
observations:

\noindent \textsl{Observation 1.} Set $K_{i}:=x_i+1$ for
$i=1,\ldots,k-1.$ Then $K_i(F)=x_{i+1}+1$.

\noindent \textsl{Observation 2.}  If we set $L_i:=1+x_i+x_{i+1}$,
then for all $1\leq i\leq k-1$, $L_i(F)=1+x_{i+1}+x_{i+2}=L_{i+1}$,
and
$$L_{k-1}(F)=1+x_k+x_{k+1}=1+x_k+\dps{\frac{a+\sum_{i=2}^{k}
x_i}{x_1}}=\frac{a+\sum_{i=1}^k x_i+x_1x_k}{x_1}.
$$
For this reason
\begin{itemize}
  \item[(a)] If $M_1:=\prod_{i=1}^{k-1}(1+x_i+x_{i+1})$, then
  $$
M_1(F)=\frac{(a+\sum_{i=1}^{k}x_i+x_1x_k)\left(\prod_{i=2}^{k-1}(1+x_i+x_{i+1})\right)}{x_1}.
  $$
  \item[(b)] Setting $M_m:=\prod_{i=1,i\neq m-1,m}^{k-1}(1+x_i+x_{i+1})$ for $2\leq m\leq k-2$, we obtain
  $$
  M_m(F)=\left(\prod_{i=2,i\neq
  m,m+1}^{k-1}(1+x_i+x_{i+1})\right)(a+\sum_{i=1}^k
  x_i+x_1x_k)/x_1.
  $$
  \item[(c)] If $M_{k-1}:=\prod_{i=1}^{k-3}(1+x_i+x_{i+1})$, then
  $M_{k-1}(F)=\prod_{i=2}^{k-2}(1+x_i+x_{i+1})$.
\end{itemize}

 \noindent \textsl{Observation 3.} Set
$N=a+\sum\limits_{i=1}^{k-1} x_i-x_2x_k$, then

\begin{equation*}
\begin{split}
N(F)&=a+\sum\limits_{i=2}^{k} x_i-x_3x_{k+1}=a+\sum\limits_{i=2}^{k}
x_i-x_3\dps {\frac{a+\sum_{i=2}^{k} x_i}{x_1}}\\
&=\dps{\frac{a+\sum_{i=2}^{k} x_i}{x_1}}\left(x_1-x_3\right)=x_{k+1}(x_1-x_3).
\end{split}
\end{equation*}

\noindent \textsl{Observation 4.} Set $R=a+\sum\limits_{i=1}^{k}
x_i+x_1x_k$, then
$$
R(F)=a+\sum\limits_{i=2}^{k}
x_i+x_{k+1}+x_2x_{k+1}=\dps{\frac{a+\sum_{i=2}^{k}
x_i}{x_1}}\left(1+x_1+x_2\right)=x_{k+1}(1+x_1+x_2).
$$

\noindent \textsl{Observation 5.} (a) For all $2\leq i\leq k-2$ set
$S_i=x_{i-1}-x_{i+1}$, then $S_i(F)=x_i-x_{i+2}$.

(b) Set $S_{k-1}=x_{k-2}-x_{k}$, then
$$S_{k-1}(F)=x_{k-1}-x_{k+1}=x_{k-1}-\dps{\frac{a+\sum_{i=2}^{k}
x_i}{x_1}}=-\dps{\frac{a+\sum_{i=2}^{k}x_i-x_1x_{k-1}}{x_1}}.$$

If we now consider the seed $X_1$ given by equation (\ref{X1}),
using Observations 1, 2a, and 3 we obtain that
$$X_2=\dps{\frac{(x_2+1)\left[\prod_{i=3}^{k-1} (1+x_i+x_{i+1})\right](a+\sum_{i=1}^{k}
x_i+x_1x_k)(x_1-x_3)}{\prod_{i=1,i\neq 2}^k x_i}}.$$

Now applying systematically Observations 1, 2b, 4 and 5a, we obtain
that for $2\leq m\leq k-1$ the component $X_m=X_{m-1}(F)$ is given
by equation (\ref{Xm}).

Observe that in particular
$$
X_{k-1}=X_{k-2}(F)=\dps{\frac{(x_{k-1}+1)\left[\prod_{i=1}^{k-3}
(1+x_i+x_{i+1})\right](a+\sum_{i=1}^{k}
x_i+x_1x_k)(x_{k-2}-x_{k})}{\prod_{i=1}^{k-2} x_i}},
$$hence the term $L_{k-1}=1+x_{k-1}+x_k$ does not appear and so,
in order to compute $X_k$ we need to use Observations 1, 2c, 4 and
5b, obtaining the expression of $X_k$ given by (\ref{Xk}).

\vspace{0.4cm}

\noindent \textbf{Second step (compatibility condition
(\ref{comp})).} A simple computations shows that
$$
X_k(F)=\frac{\mathbf{A}}{x_1\left(\prod_{i=1}^k x_i\right)},
$$
where $$\mathbf{A}=-\left(a+\sum_{i=1}^{k}
x_i\right)\left[\prod_{i=2}^{k-1}L_i\right] \left(a+\sum_{i=2}^k
x_i+x_1\left(a+\sum_{i=3}^k x_i -x_2x_k\right)\right).$$

Another computation gives that
$$
-\dps{\left(\frac{a+\sum_{i=2}^{k}x_i}{x_1^2}\right)}X_1+\dps{\frac{1}{x_1}\left[\sum_{i=2}^k
X_i \right]}=\frac{\mathbf{B}}{x_1\left(\prod_{i=1}^k x_i\right)},
$$
where $$\begin{array}{rl}
          \mathbf{B}= & -(x_1+1)\left(a+\sum_{i=2}^{k}
x_i\right)\left[\prod_{i=2}^{k-1}L_i\right] \left(a+\sum_{i=1}^{k-1}
x_i-x_2x_k\right)+\\ & \left(a+\sum_{i=1}^k x_i +x_1x_k\right)
\mathbf{C}
            -x_k(x_k+1)\left(a+\sum_{i=2}^k x_i
-x_1x_{k-1}\right)\left[\prod_{i=1}^{k-2}L_i\right],
        \end{array}
$$
and
$$
\mathbf{C}=\sum_{m=2}^{k-1} x_m(x_m+1)S_m M_m.
$$
Recall that $L_i=1+x_i+x_{i+1}$ for all $i=1,\ldots,k-1$,
$S_m=x_{m-1}-x_{m+1}$ and $M_m=\prod_{i=1,i\neq
m-1,m}^{k-1}(1+x_i+x_{i+1})$ for $2\leq m\leq k-1$. Therefore we
want to prove that $\mathbf{A}=\mathbf{B}$.

\noindent  \textsl{Step 2a.} First we show that $\mathbf{B}$
contains $L_2$ and $L_3$ as a factors. To see this, it suffices to
check that $L_2$ and $L_3$ are factors of $\mathbf{C}$. Observe that
if $L_2=0$, then we can write $L_3=x_4-x_2$, and $x_3=-(1+x_2)$
hence $x_1-x_3=x_1+x_2+1$. Therefore
$$\begin{array}{rl}
    \left.\mathbf{C}\right|_{\{L_2=0\}}=& \sum_{m=2}^{3} x_m(x_m+1)S_m M_m\\
    =&
\sum_{m=2}^{3} x_m(x_m+1)(x_{m-1}-x_{m+1}) \left(\prod_{i=1,i\neq
m-1,m}^{k-1}(1+x_i+x_{i+1})\right)\\
    =& x_2(1+x_2)(x_1-x_3)(x_4-x_2)\left(\prod_{i=4}^{k-1}(1+x_i+x_{i+1})\right) \\
    & +x_2(1+x_2)(x_2-x_4)(1+x_1+x_2)\left(\prod_{i=4}^{k-1}(1+x_i+x_{i+1})\right)=0.
  \end{array}
$$
So it is a factor of $\mathbf{C}$.

If $L_3=0$, then $L_2=x_2-x_4$, and $x_3-x_5=-(1+x_4+x_5)$. Hence
$$\begin{array}{rl}
    \left.\mathbf{C}\right|_{\{L_3=0\}}=& \sum_{m=3}^{4} x_m(x_m+1)S_m
    M_m=\\ =&
\sum_{m=3}^{4} x_m(x_m+1)(x_{m-1}-x_{m+1}) \left(\prod_{i=1,i\neq
m-1,m}^{k-1}(1+x_i+x_{i+1})\right)\\
    =& x_4(1+x_4)(x_2-x_4)(1+x_1+x_2)\left(\prod_{i=4}^{k-1}(1+x_i+x_{i+1})\right) \\
    & -x_4(1+x_4)(1+x_4+x_5)(1+x_2+x_3)(1+x_1+x_2)\times\\
    &\times \left(\prod_{i=5}^{k-1}(1+x_i+x_{i+1})\right)=0.
  \end{array}
$$
Hence $L_3$ is a factor of $\mathbf{C}$.

The above results imply that $$ \mathbf{C}=L_2 L_3
Q_k(x_1,x_2,x_4,\ldots,x_k),$$ (observe that $x_3$ does not appear
in the expression of $Q_k$). Hence $L_2$ and $L_3$ are factors in
the expression of $\mathbf{B}$, and then

$$\begin{array}{rl}
          \mathbf{B}= & L_2L_3\left[-(x_1+1)\left(a+\sum_{i=2}^{k}
x_i\right)\left[\prod_{i=4}^{k-1}L_i\right] \left(a+\sum_{i=1}^{k-1}
x_i-x_2x_k\right)\right.\\
           & \left.+\left(a+\sum_{i=1}^k x_i +x_1x_k\right) Q_k-x_k(x_k+1)\left(a+\sum_{i=2}^k x_i
-x_1x_{k-1}\right)L_1\left[\prod_{i=4}^{k-2}L_i\right]\right].
        \end{array}
$$
Furthermore
$$\begin{array}{rl}
Q_k=&-x_2(1+x_2)\left(\prod_{i=4}^{k-1}L_i\right)+(x_2-x_4)L_1
\left(\prod_{i=4}^{k-1}L_i\right)+x_4(1+x_4)L_1
\left(\prod_{i=5}^{k-1}L_i\right)+\\
&+\sum_{m=5}^{k-1}x_m(x_m+1)(x_{m-1}-x_{m+1})L_1
\left(\prod_{i=4,j\neq m-1,m}^{k-1}L_j\right).
  \end{array}
$$

\noindent  \textsl{Step 2b.} Now we state the following claim, that
will we proved at the end of the proof.

\vspace{0.4cm}

 \noindent \textbf{Claim:} For $k\ge6,$ $Q_k(x_1,x_2,x_4,\ldots,x_k)=\left(\prod_{i=4}^{k-2}
L_i\right)\left[x_1x_2L_{k-1}-x_{k-1}x_{k}L_1\right]$.

\vspace{0.2cm}

By using the claim, $\mathbf{B}=
\left[\prod_{i=2}^{k-2}L_i\right]\mathbf{D}$, where
$$\begin{array}{rl}
          \mathbf{D}= & -(x_1+1)\left(a+\sum_{i=2}^{k}
x_i\right) \left(a+\sum_{i=1}^{k-1}
x_i-x_2x_k\right)L_{k-1}\\
&+\left(a+\sum_{i=1}^k x_i +x_1x_k\right) \left[x_1x_2L_{k-1}-x_{k-1}x_{k}L_1\right]\\
           & -x_k(x_k+1)\left(a+\sum_{i=2}^k x_i
-x_1x_{k-1}\right)L_1.
        \end{array}
$$

\noindent  \textsl{Step 2c.} Observe that
$$\begin{array}{rl}
&-x_k(x_k+1)\left(a+\sum_{i=2}^k x_i
-x_1x_{k-1}\right)L_1-x_{k-1}x_{k}L_1 \left(a+\sum_{i=1}^k x_i
+x_1x_k\right)\\
=& -L_1 x_k \left[(x_k+1)\left(a+\sum_{i=2}^k x_i
-x_1x_{k-1}\right)-x_{k-1}\left(a+\sum_{i=1}^k x_i
+x_1x_k\right)\right]\\
=&-L_1 x_k \left[x_{k-1}\left(a+\sum_{i=1}^k x_i +x_1x_k-x_1-x_1
x_{k}\right)+(1+x_k)\left(a+\sum_{i=2}^k x_i\right) \right]\\
=&-L_1 x_k \left[x_{k-1}\left(a+\sum_{i=2}^k x_i
\right)+(1+x_k)\left(a+\sum_{i=2}^k x_i\right) \right]\\
=&-L_1 x_k L_{k-1}\left(a+\sum_{i=2}^k x_i\right).
\end{array}
$$
Thus $\mathbf{D}$ also contains the factor $L_{k-1}$, and therefore
$\mathbf{B}=\left[\prod_{i=2}^{k-1}L_i\right] \mathbf{E}$, where
$$\begin{array}{rl}
\mathbf{E}=& \left[-(x_1+1)\left(a+\sum_{i=2}^{k} x_i\right)
\left(a+\sum_{i=1}^{k-1}
x_i-x_2x_k\right)\right.\\
&\left.+\left(a+\sum_{i=1}^k x_i +x_1x_k\right) x_1x_2-x_k L_1
\left(a+\sum_{i=2}^k x_i\right)\right].
\end{array}
$$

\noindent  \textsl{Step 2d.} Is not difficult to check that
$\mathbf{E}$ is a quadratic polynomial in $x_3$. Now we see that
$\mathbf{E}$ vanishes either when  $x_3=x_{3,1}:=-a-\sum_{i=1,i\neq
3}^{k}x_i$ (that is when $a+\sum_{i=1}^k x_i=0$), and when
$$
x_3=x_{3,2}:=-\left(a+\sum_{i=4}^{k-1}x_i+\frac{x_2(1-x_1x_k)}{1+x_1}\right).
$$
In the first case we have
$$\begin{array}{rl}
\left.\mathbf{E}\right |_{\{a+\sum_{i=1}^k x_i=0\}}=&x_1(x_1+1)(-x_k-x_2x_k)+x_1^2x_2x_k+L_1x_1x_k \\
=&x_1x_k\left[-(x_1+1)(x_2+1)+x_1x_2+1+x_1+x_2\right]\\
=&x_1x_k\left[-1-x_1-x_2-x_1x_2+x_1x_2+1+x_1+x_2\right]=0.
\end{array}
$$

In the second case, since $a+\sum_{i=3}^{k-1}
x_i=x_2(x_1x_k-1)/(x_1+1)$, we have
$$
\begin{array}{rl}
\left. \mathbf{E}\right
|_{\{x_3=x_{3,2}\}}=&-\left(x_2+\dps{\frac{x_2(x_1x_k-1)}{x_1+1}}\right)(x_1+1)\left(
\dps{\frac{x_2(x_1x_k-1)}{x_1+1}}+x_1+x_2-x_k-x_2x_k\right)\\
&+x_1x_2\left(x_1+x_2+\dps{\frac{x_2(x_1x_k-1)}{x_1+1}}+x_1x_k\right)\\
&-L_1
x_k\left(x_2+\dps{\frac{x_2(x_1x_k-1)}{x_1+1}}\right)=:\dps{\frac{1}{x_1+1}}\,\mathbf{F},
\end{array}
$$
where (after some computations) $$
\begin{array}{rl}
\mathbf{F}=&-x_2x_1(x_k+1)L_1(x_1-x_k)+x_1^2x_2L_1(x_k+1)-L_1x_1x_2x_k(x_k+1)\\
=&x_1x_2L_1(x_k+1)[-(x_1-x_k)+x_1-x_k]=0.
\end{array}
$$
Therefore $\left. \mathbf{E}\right |_{\{x_3=x_{3,2}\}}=0$.

In summary,  as a quadratic polynomial in $x_3$, $\mathbf{E}$
factorizes as
$$\begin{array}{rl}
\mathbf{E}=&-(x_1+1)\left(x_3+a+\sum_{i=1,i\neq
3}^{k}x_i\right)\left(x_3+a+\sum_{i=4}^{k-1}x_i+(x_2(1-x_1x_k))/(x_1+1)\right)\\
=&-\left(a+\sum_{i=1}^k x_i\right) \left(\left(
x_3+a+\sum_{i=4}^{k-1} x_i\right)(x_1+1)+x_2(1-x_1x_k)\right)
\\
=&-\left(a+\sum_{i=1}^k x_i\right) \left((x_1+1)\left(
a+\sum_{i=3}^{k-1} x_i\right)+x_2-x_1x_2x_k\right)\\
=&-\left(a+\sum_{i=1}^k x_i\right) \left(a+\sum_{i=2}^{k-1}x_i+x_1
\left( a+\sum_{i=3}^{k-1} x_i-x_2x_k\right)\right).
\end{array}
$$

So, finally, we get that
$$\begin{array}{rl}
\mathbf{B}=&-\left[\prod_{i=2}^{k-1}L_i\right]\left(a+\sum_{i=1}^k
x_i\right) \left(a+\sum_{i=2}^{k-1}x_i+x_1 \left( a+\sum_{i=3}^{k-1}
x_i-x_2x_k\right)\right)=\mathbf{A},
\end{array}
$$
as we wanted to show.

To end the proof it only remains to prove the claim. We proceed by
induction. That it is true when $k=6$ is straightforward.
 Assume now that the claim is true for
$Q_k$, then
$$\begin{array}{rl}
    Q_{k+1}= & -x_2(1+x_2)\left(\prod_{i=4}^{k}L_i\right)+(x_2-x_4)L_1\left(\prod_{i=4}^{k}L_i\right)
+x_4(1+x_4)L_1\left(\prod_{i=5}^{k}L_i\right)+\\
     & +\sum_{m=5}^{k-1}x_m(1+x_m)(x_{m-1}-x_{m+1})L_1
\left(\prod_{i=4,i\neq m-1,m}^{k}L_i\right) \\
    = & L_k\,
    Q_k+L_1\left(\prod_{i=4}^{k-2}L_i\right)\left(x_k(1+x_k)(x_{k-1}-x_{k+1})\right).
  \end{array}
$$
By using the hypothesis of induction, we obtain that
$$\begin{array}{rl}
    Q_{k+1}= & L_k\left(\prod_{i=4}^{k-2}L_i\right)\left[x_1 x_2 L_{k-1}-x_{k-1}x_k
    L_1\right]+\\
    & L_1\left(\prod_{i=4}^{k-2}L_i\right)\left(x_k(1+x_k)(x_{k-1}-x_{k+1})\right)=\\
    & \left(\prod_{i=4}^{k-2}L_i\right)\left(x_1x_2 L_k
    L_{k-1}+L_1\left[x_k(1+x_k)(x_{k-1}-x_{k+1})-x_{k-1}x_k
    L_k\right]\right).
  \end{array}
$$
Thus
$$\begin{array}{rl}
    Q_{k+1}=& \!\!\!\! L_k\left(\prod_{i=4}^{k-2}L_i\right)\left[x_1 x_2 L_{k-1}-x_{k-1}x_k
    L_1\right]+ L_1\left(\prod_{i=4}^{k-2}L_i\right)\left(x_k(1+x_k)(x_{k-1}-x_{k+1})\right)\\
    =&\!\!\!\! \left(\prod_{i=4}^{k-2}L_i\right)\left(x_1x_2 L_k
    L_{k-1}+L_1\left[x_k(1+x_k)(x_{k-1}-x_{k+1})-x_{k-1}x_k
    L_k\right]\right).
  \end{array}
$$

An easy computation shows that
$x_k(1+x_k)(x_{k-1}-x_{k+1})-x_{k-1}x_k L_k=-L_{k-1}x_k x_{k+1}$,
and the result follows. Therefore, Theorem \ref{propoSL} is proved.

\section{Geometrical issues in the odd case.}\label{geoiss}

Before proving Theorem \ref{propnova} and Propositions \ref{propoL4}
and \ref{propoL5}  we need some preliminary results. Recall that a
map $H$ which is a first integral for
$G^p:=G\circ\overset{_{_{\,\,p)}}}{\cdots} \circ G$ is also called
sometimes {\it a $p-$first integral}, or simply, for short, a {\it
$p-$integral}, of $G,$ see \cite{RQ}.

The following lemma, which is very easy to prove,  gives light on
one utility of $p-$first integrals, specially if they are not
symmetric functions of their arguments.
\begin{lem}\label{metode}
Let  $H$ be a $p$-integral of a map $G.$ Then  then for any
symmetric function of $p$ variables $S$, the function
$V_S:=S(H,H\circ G,H\circ G^2,\ldots, H\circ G^{p-1})$ is a first
integral of $G$.
\end{lem}

\begin{lem}\label{lemag} Set $k=2\ell+1, $ ${\bf x}=(x_1,x_2,\ldots,x_k)\in {\mathcal Q}^+$
and let $W$ be the function given in~(\ref{W}). Then, if we define
the polynomial
\begin{equation*}\label{eq1}
\begin{array}{rl}
 Z({\bf x}):&= (\prod_{i=1}^{k} x_i) [W({\bf x})-W(F({\bf x}))]=\\
   &=\prod_{j=0}^{\ell}x_{2j+1}(x_{2j+1}+1)-
\left(a+\sum_{i=1}^{2\ell
+1}x_i\right)\prod_{j=1}^{\ell}x_{2j}(x_{2j}+1),
\end{array}
\end{equation*} it holds that $Z(F({\bf x}))=\det(DF({\bf x})) Z({\bf x}).$
\end{lem}

\dem In Theorem \ref{propnova} (a) we will prove  that $W$ is a
2-integral of $F.$ Thus, if we define $\tilde{Z}:=W-W(F),$ then
\begin{equation}\label{eqgtilde0}
\tilde{Z}(F({\bf x}))=-\tilde{Z}({\bf x}). \end{equation} From the
above equality it is clear that $\{{\bf x}\in{\cal Q}^+ \,:\,
\tilde{Z}({\bf x})=0\}=\{{\bf x}\in{\cal Q}^+ \,:\, Z({\bf x})=0\}$
is an invariant hypersurface by $F.$ Notice also that if we define $
\Pi({\bf x}):=\prod_{j=1}^{k}x_{j}, $ it holds that
\begin{equation}\label{eqgtilde1}
\Pi(F({\bf x}))=\frac{a+\sum_{i=2}^k x_i}{x_1^2}\Pi({\bf x})
=-\det(DF({\bf x}))\, \Pi({\bf x}).
\end{equation}
Since $Z({\bf x})=\Pi({\bf x})\tilde{Z}({\bf x}),$ by using
equalities (\ref{eqgtilde0}) and (\ref{eqgtilde1}), we obtain that
$$
Z(F({\bf x}))=\Pi(F({\bf x}))\tilde{Z}(F({\bf x}))=\det(DF({\bf x}))\Pi({\bf x})\tilde{Z}({\bf x})
=\det(DF({\bf x})){Z}({\bf x}),
$$
as we wanted to see.\qed

\noindent{\bf Proof of Theorem \ref{propnova}.} (a) The proof of the
equality  $ W(F^2({\bf x}))=W({\bf x})$ is straightforward. We have
obtained the expression (\ref{W}) inspired in the results of
\cite{MT}. In that paper it is proved that
$$\frac{(y_{n}+1)(y_{n+2}+1)}{y_{n+1}}=\frac{(y_{n+2}+1)(y_{n+4}+1)}{y_{n+3}},$$
where $\{y_n\}$ is the sequence given by the 3-rd order Lyness'
recurrence $y_{n+3}=(a+y_{n+1}+y_{n+2})/y_{n}.$ Notice that this
property is equivalent to say that for $k=3,$ $W$ is a  2-integral
for~$F.$

(b-c) By applying Lemma \ref{metode} with $S(u,v)=u+v$ and
$S(u,v)=uv,$ we obtain the first integrals $V_3$ and $V_1,$
respectively. The functionally independence of $V_1,V_2$ and $V_3,$
for $\ell\ge2,$ follows from straightforward computations and it is
already established in \cite{G}.

(d-f) From Lemma \ref{lemag} we know that $Z(F({\bf
x}))=\det(DF({\bf x})){Z}({\bf x}),$ where recall that $Z({\bf x})=
(\prod_{i=1}^{k} x_i) [W({\bf x})-W(F({\bf x}))].$  Note also that
$$\det(DF({\bf
x}))=(-1)^k\displaystyle{\frac{a+x_2+\cdots+x_{k-1}}{x_k^2}}.$$
Since when $k$ is odd $\det(DF)<0$ on $\mathcal{Q}^+$, equation
$Z(F)=\det(DF) Z$ means that
$$\mathcal{G}=\{{\bf x}\in{\cal Q}^+ \,:\,
W({\bf x})=W(F({\bf x}))\}=\{{\bf x}\in{\cal Q}^+ \,:\,
Z({\bf x})=0\}
$$ is invariant by $F$ and that $F$
maps the region $\{{\bf x}\,:\,Z({\bf x})>0\}$ into the region
$\{{\bf x}\,:\,Z({\bf x})<0\}$ and viceversa. Furthermore it implies
that the dynamics of $F^2$ on each of these sets are conjugated,
being the map $F$ itself the conjugation. Moreover, any periodic
orbit with odd period must lie in $\sg$, as we wanted to see.

(g-h) By using the Change of Variables Theorem it is easy to see
that if $G$ is a diffeomorphism of $\mathcal{U}$, and on this region
$\mu$ is a positive function that satisfies $\mu(G({\bf
x}))=\det(DG({\bf x}))\mu ({\bf x}),$ then
$$m(B)=\int_{B} \frac{1}{\mu({\bf x})}\,d{\bf x}$$
is an invariant measure for $G.$ By Lemma \ref{lemag} we know that
$Z(F({\bf x}))=\det(DF({\bf x}))Z({\bf x})$ and by  equality
(\ref{eqgtilde1}), that $\Pi(F({\bf x}))=-\det(DF({\bf x}))\Pi({\bf
x}).$ By using these results we have that $Z(F^2({\bf
x}))=\det(DF^2({\bf x}))Z({\bf x})$ and $\Pi(F^2({\bf
x}))=\det(DF^2({\bf x}))\Pi({\bf x}),$ being both equalities in the
corresponding domains, which are invariant by $F^2$. Hence (f) and
(g) follow. \qed

\subsection{Order reduction}\label{dinreduida}

Recall that in the previous section we have seen that when
$k=2\ell+1$, the regions $\{{\bf x}\,:\,Z({\bf x})>0\}$ and $\{{\bf
x}\,:\,Z({\bf x})<0\}$ are invariant by $F^2$, and the dynamics on
both region are conjugated. This observation allows us to give a new
application of the invariant $W.$ We can  reduce the study of the
dynamics of $F$ on $\{{\bf x}\,:\,Z({\bf x})\ne0\}$ to the study of
a new $(k-1)$-dimensional map, having one more parameter.

For instance for $n=3$, we get that $W(x,y,z)=(x+1)(z+1)/y$, and
hence any admissible level surface $\{{\bf x}\,:\,W({\bf x}) =w\}$,
$w\ne0,$ can be described as $y=k(x+1)(z+1)$, where $k=1/w$.
Therefore
 $$\left. F^2\right
|_{\{W=w\}}(x,y,z)=\left(z,\frac{a+z+k(x+1)(z+1)}{x},\frac{a+k+z(k+1)}{k
x(z+1)}\right),
$$
and we can reduce the study of the dynamics of $F^2$ to the study of
the reduced map
\begin{equation*}\label{FG}\tilde{F}_2(x,z)=\left(z,\frac{a+k+z(k+1)}{k
x(z+1)}\right),
\end{equation*}
or, equivalently, the study of the second order difference equation
\begin{equation}\label{genlyn}
y_{n+2}=\frac{a+k+y_{n+1}(k+1)}{ky_{n}(y_{n+1}+1)},
\end{equation}
as in \cite[Eqs. (8) and (9)]{MT}. The dynamics of this difference
equation is studied in \cite[Ex. 3]{CGM07-2}. Equation
(\ref{genlyn}) is  sometimes called {\it generalized Lyness'
recurrence}, see \cite{K}.

For $n=5$, the integral of $F^2$ is
$W(x,y,z,t,s)=(x+1)(z+1)(s+1)/(yt)$. Again any admissible level
surface $\{{\bf x}\,:\,W({\bf x}) =w\}$, $w\ne0,$ can be described
as $t=k(x+1)(z+1)(s+1)/y$, where $k=1/w$. Therefore proceeding as
before we can reduce the study of the dynamics of $F^2$ to the study
of the reduced map
$$
\tilde{F}_2(x,y,z,s)=\left(y,z,t,{\frac {p_2(z,s;k)
{x}^{2}+p_1(y,z,s;a,k)
 x+p_0(y,z,s;a,k) }{x{y}^{2}}}
\right)
$$
where $p_2(z,s;k)=k \left( s+1 \right) \left( z+1 \right)$,
$p_1(y,z,s;a,k)= 2\,k
 \left( s+1 \right)  \left( z+1 \right) +y \left( a+s+z \right)
$, and $p_0(y,z,s;a,k)=k \left( s+1 \right)  \left( z+1 \right) +y
\left( z+s+a+y
 \right)$,   or equivalently to the difference equation
$$
y_{n+4}=\frac {p_2(y_{n+2},y_{n+3};k)
{y_{n}}^{2}+p_1(y_{n+1},y_{n+2},y_{n+3};a,k)
 y_{n}+p_0(y_{n+1},y_{n+2},y_{n+3};a,k) }{y_{n}{y_{n+1}}^{2}}.
$$
Clearly, the described procedure can be generalized to higher
dimensions.

\section{Dynamics of the low dimensional cases}\label{lowlow}

This subsection is devoted to prove Propositions \ref{propoL4} and
\ref{propoL5}.

Along the section we will use  a straightforward consequence of the
following result:
\begin{teo}[Bastien and Rogalsky, \cite{BR}] Let $\bar\bx$ be the fixed point of $F$ in ${\mathcal Q}^+.$
For any $h>V_1(\bar\bx)$, the level sets  $\{V_1=h\}\cap{\mathcal
Q}^+$ are homeomorphic to $\mathbb{S}^{k-1}$.
\end{teo}
\begin{corol}\label{compa}
Let $K\ne\emptyset$ be  by the intersection of some level sets of
different first integrals of $F,$ including $V_1$ among them. Then
$K\cap{\mathcal Q}^+$ is a  compact  set, invariant by $F.$
\end{corol}

For the sake of completeness and to compare with the cases $k=4,5$,
we start by recalling some results for the case $k=3,$ most of them
already proved in \cite{CGM07-1}.

\subsection{The $3$-dimensional map}\label{cas3}

For $k=3,$
\begin{equation*}\label{lyn}
F(x,y,z)=\left(y,z,\dps{\frac{a+y+z}{x}}\right),
\end{equation*}
the  Lie symmetry given in Theorem~\ref{propoSL}, is
$$
\begin{array}{rl}
{\bf X}_{3} & :=\left[x(x+1)(1+y+z)(a+x+y-yz)\,\frac{\partial}{\partial x}+y(y+1)(x-z)(a+x+y+z+xz)\,
\frac{\partial}{\partial y}\right.\\
&\left.+z(z+1)(1+x+y)(a+y+z-xy)\,\frac{\partial}{\partial
z}\right]/(xyz)
\end{array}
$$
and since  ${\bf X}_{3}(V_i)=0$, for $i=1,2,$ the functions $V_1$
and $V_2$ given in (\ref{v1}) and (\ref{v2}) are first integrals for
${\bf X}_3$ and $F.$ Also, by Theorem \ref{propnova}, $ W(x,y,z):=
\left( x+1 \right) \left( z+1 \right)/y $ is a 2-integral of $F$; $
\sg=\{{\bf x}\in{\mathcal Q}^+\,:\, Z({\bf x})=0\}$ is invariant by
$F,$ where
$$ Z({\bf x})= x(x+1)z(z+1)-(a+x+y+z)y(y+1);$$
and $F$ maps $\sg^+:=\{{\bf x}\in{\mathcal Q}^+\,:\, Z({\bf x})>0\}$
into $\sg^-:=\{{\bf x}\in{\mathcal Q}^+\,:\, Z({\bf x})<0\},$ and
viceversa.

Let $\bar\bx$ be the fixed point in ${\mathcal Q}^+$, of $F$. Set
$h>V_1(\bar\bx)$, $k>V_2(\bar\bx)$. In \cite{CGM07-1} it is proved
that $\{V_1=h\}\cap {\mathcal Q}^+$ and $\{V_2=k\}\cap {\mathcal
Q}^+$ are diffeomorphic to two dimensional spheres, and  their
transversal intersections are formed by exactly two disjoints
curves, both diffeomorphic to circles. Their non transversal
intersections correspond to:

\noindent (a) The $2$--periodic points of $F$ (which are equilibrium
points of $\mathbf{X}_3$) given by the curve ${\cal
L}:=\{(x,(x+a)/(x-1),x)\,|\,x>1\}.$

\noindent (b) The levels placed at $\mathcal{G}$. Those placed at
$\mathcal G\setminus\{\bar\bx\}$ are formed by exactly one curve,
diffeomorphic to a circle.

Finally, since ${\bf X}_{3}$ is also a  Lie symmetry of $F^2$, as a
consequence of \cite[Thm. 18]{CGM07-1} or \cite[Thm. 1]{CGM07-2}, we
know that the map $F^2$ restricted to each of the sets
$\{V_1=h\}\cap\{V_2=k\},$ which is not simply a fixed point of
$F^2,$ is conjugated to a rotation. Further discussion about the
possible rotation numbers can be found in \cite{CGM07-1}.

\subsection{The $4$-dimensional map.}

\noindent {\bf Proof of Proposition \ref{propoL4}}. Equations
(\ref{X1})--(\ref{Xk}) give the following Lie symmetry for the $4$-dimensional
Lyness' map:

\begin{equation*}\label{SL4}
\begin{array}{rl}
 {\bf X}_{4}=&\left[{ { x\left( x+1 \right)  \left( 1+y+z \right) \left(
1+z+t
 \right)  \left( a+x+y+z-yt \right) }} {\frac{\partial}{\partial x}}\right.\\
 &+{{  y\left( y+1 \right)\left( 1+z+t
 \right) \left( a+x+y+z+t +xt\right)
   \left( x-z \right)}}{\frac{\partial}{\partial y}}\\
   & +{ { z\left( z+1 \right)  \left( 1+x+y \right)\left( a+x+y+z+t +xt
\right)
  \left( y-t \right) }}{\frac{\partial}{\partial z}} \\
&\left. -{ { t\left( t+1 \right) \left(1 +x+y\right)\left( 1+y+z
\right)\left( a+y+z+t -xz\right)
   }}{\frac{\partial}{\partial
t}}\right]/(xyzt). \end{array}\end{equation*}

A straightforward computation shows that the Lie symmetry satisfies
${\bf X}_{4}(V_i)=0$, for $i=1,2$. Hence,  the orbits of both $F$
and ${\bf X}_{4}$ generically lie in a two dimensional surface of
the form  $I_{h,k}:=\{V_1=h\}\cap\{V_2=k\}\cap {\cal Q}^+.$

Let $C_{h,k}$  be a connected component of  $I_{h,k}$.  From
Corollary~\ref{compa}  we know that $C_{h,k}$ is compact. If
$\{V_1=h\}$ and $\{V_2=k\}$ intersect transversally on $C_{h,k}$,
then for all points
 $$ {\rm Rank}\left(\begin{array}{cccc}
                  (V_1)_x & (V_1)_y & (V_1)_z & (V_1)_t \\
                  (V_2)_x & (V_2)_y & (V_2)_z & (V_2)_t
                \end{array}
\right)=2.
$$ This fact implies that the dual $2$--form associated to the $2$--field
$$
\begin{array}{rl}
{\nabla} V_1 \wedge {\nabla} V_2=&[(V_1)_x(V_2)_y
-(V_1)_y(V_2)_x]\frac{\partial}{\partial x}\wedge
\frac{\partial}{\partial
y}+[(V_1)_x(V_2)_z-(V_1)_z(V_2)_x]\frac{\partial}{\partial x}\wedge
\frac{\partial}{\partial z}\\
&+[(V_1)_x(V_2)_t-(V_1)_t(V_2)_x]\frac{\partial}{\partial x}\wedge
\frac{\partial}{\partial
t}+[(V_1)_y(V_2)_z-(V_1)_z(V_2)_y]\frac{\partial}{\partial y}\wedge
\frac{\partial}{\partial z}\\
&+[(V_1)_y(V_2)_t-(V_1)_t(V_2)_y]\frac{\partial}{\partial y}\wedge
\frac{\partial}{\partial t}+
[(V_1)_z(V_2)_t-(V_1)_t(V_2)_z]\frac{\partial}{\partial z}\wedge
\frac{\partial}{\partial t},
\end{array}
$$ given by
$$\begin{array}{rl}
\omega=&
[(V_1)_z(V_2)_t-(V_1)_t(V_2)_z]dxdy+[(V_1)_z(V_2)_x-(V_1)_x(V_2)_z]dxdz\\
&+[(V_1)_z(V_2)_y-(V_1)_y(V_2)_z]dxdt+
[(V_1)_t(V_2)_x-(V_1)_x(V_2)_t]dydz\\
&+[(V_1)_t(V_2)_y-(V_1)_y(V_2)_t]dydt+[(V_1)_x(V_2)_y
-(V_1)_y(V_2)_x]dzdt
\end{array}
$$
is nonzero at every point of  $C_{h,k}$, and therefore it is
orientable, see \cite[Sec. 2.5]{AM} or \cite[Sec. X, \S 1.]{L}.

Let ${\bf X}_{4}^{h,k}$ denote the restriction of the vector field
${\bf X}_{4}$ to the invariant surface $C_{h,k}$. Some computations
show that the unique equilibrium point of ${\bf X}_4$ in ${\cal
Q}^+$ is the fixed point of $F.$ Hence ${\bf X}_{4}^{h,k}$ has no
equilibrium points, and therefore the Poincar\'e--Hopf formula
gives:
$$
0=i({\bf X}_{4}^{h,k})=\chi(C_{h,k})=2-2g
$$
where   $i({\bf X}_{4}^{h,k})$ denotes the sum of the indices of the
equilibrium points of ${\bf X}_{4}^{h,k}$ in $C_{h,k}$ and
$\chi(C_{h,k})$ and $g$ are the Euler characteristic and the genus
of the surface $C_{h,k},$ respectively.  Hence the genus of
$C_{h,k}$ is one. An orientable, compact, connected surface of genus
one is a torus, as we wanted to prove.\qed

\subsection{The $5$-dimensional map}\label{cas5}

\noindent {\bf Proof of Proposition \ref{propoL5}.} By Theorem
\ref{propnova} we know that
$$ W(x,y,z,t,s):={\frac { \left( x+1 \right) \left( z+1 \right)
\left( s+1 \right) }{ yt}},
$$ is a 2-integral of $F.$    Moreover, the 3 functionally independent first integrals
of $F$ given in (\ref{v1}), (\ref{v2}) and (\ref{v3}), are
$$\begin{array}{l}
V_1(x,y,z,t,s)=\dps{\frac { \left( a+x+y+z+t+s \right) \left( x+1
\right)  \left( y+1
 \right)  \left( z+1 \right)  \left( t+1 \right)  \left( s+1 \right) }
{xyzts}}, \\
 V_2(x,y,z,t,s)=\dps{\frac { \left( a+x+y+z+t+s+xs \right)  \left(
1+x+y \right)  \left( 1 +y+z \right)  \left( 1+z+t \right)  \left(
1+t+s \right) }{xyzts}},\\
 V_3(x,y,z,t,s)=\dps{{\frac {  x\left( x+1 \right)z \left( z+1 \right) s \left( s+1
 \right) + \left( a+x+y+z+t+s \right)  y\left( y+1 \right) t \left( t+1
 \right)}{xyzts}}}.
\end{array}
$$

Some tedious computations show that the hypersurface $\sg$ is
 in the locus of non--trans\-ver\-sa\-lity of the three level
sets of the integrals $V_i$, $i=1,2,3$ in $Q^+$. Recall that
precisely, $\sg=\{{\bf x}\in{\mathcal Q}^+\,:\, Z({\bf x})=0\}$ is
invariant by $F,$ where
$$ Z({\bf x})= x(x+1)z(z+1)s(s+1)-(a+x+y+z+t+s)y(y+1)t(t+1),$$
and that $F$ maps $\sg^+=\{{\bf x}\in{\mathcal Q}^+\,:\, Z({\bf
x})>0\}$ into $\sg^-=\{{\bf x}\in{\mathcal Q}^+\,:\, Z({\bf
x})<0\},$ and viceversa.

 Equations
(\ref{X1})--(\ref{Xk}) give the following Lie symmetry for the $5$-dimensional
Lyness' map:

\begin{equation*}\label{SL5}
\begin{array}{rl}
{\bf X}_{5}= & \left[ x\left( x+1 \right)  \left( 1+y+z \right)  \left( 1+z+t
\right)
\left( 1+t+s \right)  \left( a+x+y+z+t-ys \right){\frac{\partial}{\partial x}} \right.\\
{}  &   +y\left( y+1 \right) \left( 1+t+s
\right)  \left( 1+z+t \right)  \left( a+x+y+z+t+s+xs \right) \left( x-z \right){\frac{\partial}{\partial y}}\\
{}     & +z\left( z+1 \right)  \left( 1+x+y \right) \left( 1+t+s
\right)
\left( a+x+y+z+t+s+xs \right)    \left( y-t \right){\frac{\partial}{\partial z}}\\
 {}      & +t \left( t+1 \right) \left( 1+x+y \right) \left(
1+y+z \right)
 \left( a+x+y+z+t+s+xs \right) \left( z-s \right) {\frac{\partial}{\partial t}}\\
{}   & \left. - s\left( s+1 \right)  \left( 1+x+y \right)  \left(
1+y+z \right)
 \left( 1+z+t \right)  \left( a+y+z+t+s -tx\right){\frac{\partial}{\partial s}} \right]/(xyzts).
\end{array}
\end{equation*}

Again,  direct computations show that ${\bf X}_{5}(V_i)=0$, for
$i=1,2,3.$ Hence the orbits of both $F$,  and ${\bf X}_{5}$ lie in a
two dimensional surface of the form
$I_{h,k,\ell}:=\{V_1=h\}\cap\{V_2=k\}\cap\{V_3=\ell\}\cap {\cal
Q}^+.$

Let $C_{h,k,\ell}$  be a connected component of $I_{h,k,\ell}$. From
Corollary~\ref{compa}  we know that $C_{h,k,\ell}$ is compact. If
$\{V_1=h\}$, $\{V_2=k\}$ and $\{V_3=k\}$ intersect transversally on
$C_{h,k,\ell}$, then for all points
 $$ {\rm Rank}\left(\begin{array}{ccccc}
                  (V_1)_x & (V_1)_y & (V_1)_z & (V_1)_t & (V_1)_s\\
                  (V_2)_x & (V_2)_y & (V_2)_z & (V_2)_t & (V_2)_s\\
                  (V_3)_x & (V_3)_y & (V_3)_z & (V_3)_t & (V_3)_s\
                \end{array}
\right)=3.
$$
Similarly than in the case $k=4,$ this fact implies that the dual
$2$--form associated to the $3$--field ${\nabla} V_1 \wedge {\nabla}
V_2\wedge {\nabla} V_3$ is nonzero at every point of $C_{h,k,\ell}$,
and therefore this set is a two-dimensional orientable manifold.

It is not difficult to check that all the equilibrium points of
${\bf X}_5$ in $\mathcal{Q}^+$ are the points of the curve
$$ {\cal L}=\left\{{\bf x}=\left(x,{\frac
{2\,x+a}{x-2}}, x, {\frac {2\,x+a}{x-2}},x\right) \quad
\mbox{with}\quad x>2 \right\},
$$ which  contains  a continuum of two periodic points  and the fixed point of
$F.$ Moreover ${\mathcal L}$ belongs to the locus of
non--transversality of the integrals $V_1$, and $V_2$ in
$\mathcal{Q}^+$.

The above observation implies that $\mathbf{X}_5$ has no equilibrium
points in any level set $I_{h,k,\ell}$ where the three first
integrals intersect transversally. Therefore  the Poincar\'e--Hopf
formula gives
$$
0=i({\bf X}_{5}^{h,k,\ell})=\chi(C_{h,k,\ell})=2-2g.
$$
for each connected component $C_{h,k,\ell}$ of such a level set
(where ${\bf X}_{5}^{h,k,\ell}$ is the restriction of ${\bf X}_{5}$
to $C_{h,k,\ell}$).

 Hence $g=1$, which implies that $C_{h,k,\ell}$ is a torus (since it is
 two-dimensional, orientable, compact, manifold of
genus one), as we wanted to proof.

Finally,  observe that the fact that $F$ maps $\sg^+$ into $\sg^-,$
and viceversa, implies that most $I_{h,k,\ell}$ have at least one
connected component on each set. In fact, it seems that each
$I_{h,k,\ell},$ not included in $\G,$ has exactly two connected
components in ${\mathcal Q}^+,$ as it happens when $k=3.$ See
Figures 2 and 3 for two illustrations of this assertion. \qed

\centerline{\includegraphics[scale=0.55]{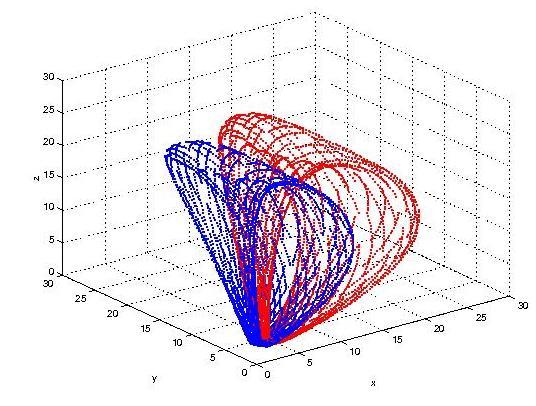}}
\begin{center}Figure 3. Projection into  $\R^3$ of the first
$10^4$ iterates of  the Lyness' map,  for $k=5$ and  $a=4,\!$
starting at  $\!(1,2,3,4,5).$\! Odd and even iterates are in
different connected components.
\end{center}

From the above proof it is clear that if some $I_{h,k,\ell}$ cuts
$\mathcal L$ then the first three integrals do not cut transversally
on it. Let us see that in general $I_{h,k,\ell}\cap
\mathcal{L}=\emptyset.$ This will be a consequence of the shape of
the function
$$
v_1(x):=V_1\left(x,{\frac {2\,x+a}{x-2}}, x, {\frac
{2\,x+a}{x-2}},x\right),\quad x>2,\, a\geq 0.
$$
This function has a global minimum at the coordinate given by the
fixed point $x=2+\sqrt{4+a}$, and $\lim\limits_{x\to 2^+}
v_1(x)=$$\lim\limits_{x\to +\infty}v_1(x)=+\infty$.

Thus, given $h> v_1(2+\sqrt{4+a})$  there are only two solutions
$x_1(h)<2+\sqrt{4+a}<x_2(h)$ of the equation $v_1(x)=h$. Now set
$$
v_i(x):=V_i\left(x,{\frac {2\,x+a}{x-2}}, x, {\frac
{2\,x+a}{x-2}},x\right),\quad x>2,\, a\geq 0,
$$
for $i=2,3.$ Then, for any value of $k$ and $\ell$ satisfying $
k\not\in\{v_2(x_1(h)),v_2(x_2(h))\}$  or $\ell\not\in
\{v_3(x_1(h)),v_3(x_2(h))\},$ we obtain that $I_{h,k,\ell}$, does
not intersect ${\cal L}$.

\section{Conclusions}\label{conclusions}

Several properties for the $k$-dimensional Lyness' map $F$ have been
given, like the existence of a  Lie symmetry for $F$ and of a new
and simple invariant for~$F^2.$  This Lie symmetry together with the
new invariant
 give  information  for $k=4$ and
$5$ about the dynamics and the  topology of the level surfaces where
the dynamics of $F$ is confined. Some general results for $k$ odd
have been also presented. However, on the contrary that happens in
the cases $k=2$ and $3$, the numerical explorations indicate that
for $k=4,5$  the orbits of the map are not contained in the orbits
of the flow of the Lie symmetry {\it with the same initial
condition}, although they are placed in the same manifold, which has
dimension $k-E(\frac{k+1}{2})$. This is an obstruction to apply the
theoretical tools developed in \cite{CGM07-2}.

Numerical simulations seem to show that for some initial conditions
the projection in $\mathbb{R}^3$ of the iterates of $F$ when $k=6,7$
fill densely a
 $2$-dimensional manifold, indicating that probably they live in
a $3$-dimensional manifold of $\mathbb{R}^6$ and $\mathbb{R}^7.$
These facts are coherent with the   conjecture of \cite{G} about the
number of independent first integrals of the Lyness' maps, and show
that  for $k\geq 6$ the dynamics are much more complicated.

When $k=2\ell,$ the simplest scenario that we imagine for the
dynamics of the $k$-dimensional Lyness' map  is that  most of the
orbits lie on invariant manifolds which are diffeomorphic to
$\ell$-dimensional tori, $S^1\times
\overset{_{_{\,\,\ell)}}}{\cdots}\times S^1.$ On the other hand,
when $k=2\ell+1,$ most of them  lie on two diffeomorphic copies of
$S^1\times \overset{_{_{\,\,\ell)}}}{\cdots}\times S^1,$ separated
by the invariant set $\mathcal G.$ Moreover these orbits jump from
one of these tori to the other one and viceversa.

In any case, much more research must  be done in order to have a
total understanding of  the dynamics and the geometrical structure
of high dimensional Lyness' maps.

\hfill

\noindent{\bf Acknowledgements.} We want to thank Guy Bastien and
Marc Rogalski for communicating their results of \cite{BR} prior to
publication. The third author is grateful to  Immaculada G\'alvez
for her kind help.

GSD-UAB and CoDALab Groups are supported by the Government of
Catalonia through the SGR program. They are also supported by DGICYT
through grants MTM2005-06098-C02-01 (first and second authors) and
DPI2005-08-668-C03-1 (third author).


\begin{thebibliography}{12}

\bibitem{AM} R. Abraham, J.E. Marsden. ``Foundations of Mechanics'' 2nd Ed.
Addison--Wesley. Redwood City, California, 1987.

\bibitem{BR1} G. Bastien and M. Rogalski.  {\sl Global behavior of the solutions of
Lyness difference equation $u_{n+2}u_n=u_{n+1}+a$}, J. Difference
Equations and Appl. 10, no. 11 (2004), 977--1003.

\bibitem{BR2} G. Bastien and M. Rogalski.  {\sl On algebraic difference equations $u_{n+2}+u_n=\psi(u_{n+1})$
in $\R$ related to a family of elliptic quartics in the plane}, J.
Math. Anal. Appl. 326 (2007), 822--844.

\bibitem{BR} G. Bastien, M. Rogalsky. {\sl Results and conjectures about global behavior
of the solutions of the order $q$ Lyness' difference equation in
$\R_{*}^+$}, work in progress. Private communication. December
13$^{{\rm th}}$, 2006.

\bibitem{BC} F. Beukers and R. Cushman.  {\sl Zeeman's monotonicity conjecture},
J. Differential  Equations 143 (1998), 191--200.

\bibitem{CGM1} A. Cima, A. Gasull, V. Ma\enya osa.
{\sl Dynamics of rational discrete dynamical systems via first
integrals}, Internat. J. Bifur. Chaos Appl. Sci. Engrg. 16 (2006),
631-645.



\bibitem{CGM07-1} A. Cima, A. Gasull, V. Ma\enya osa.
{\sl Dynamics of the third order Lyness' difference Equation}, J.
Difference Equations \& Appl. 13 (2007), 855--844.

\bibitem{CGM07-2} A. Cima, A. Gasull, V. Ma\enya osa.
{\sl Studying discrete dynamical systems through differential
equations}, J. Differential Equations 244 (2008), 630--648.




\bibitem{CM} A. Cima, F. Ma\enya osas.
{\sl Real dynamics of integrable birational maps}. In preparation,
2008.


\bibitem{BMG} L. Gardini, G.I. Bischi, C. Mira. {\sl Invariant curves
and focal points in a Lyness iterative process}, Internat. J. Bifur.
Chaos Appl. Sci. Engrg. 13 (2001), 1841--1852.

\bibitem{G} M. Gao, Y Kato, M. Ito. {\sl Some invariants for $k^{{\rm th}}$-order Lyness
equation},  Appl. Math. Lett. 17 (2004), 1183-1189.

\bibitem{GP} V. Guillemin, A. Pollack. ``Differential Topology''.
Prentice Hall. Englewood Cliffs, New Jersey, 1974.

\bibitem{GM} I. Gumovski, C. Mira. ``Recurrences and discrete
dynamic systems''. Lecture Notes in Mathematics 809. Springer
Verlag, Berlin, 1980.

\bibitem{HBQC} F.A. Haggar, G.B. Byrnes, G.R. Quispel and  H.W. Cappel.
{\sl $k$-integrals and $k$-Lie symmetries in discrete dynamical
systems,} Phys. A 233 (1996), 379--394.

\bibitem{HKY} R. Hirota, K. Kimura and  H. Yahagi. {\sl How to find the
conserved quantities of nonlinear discrete equations} J. Phys.
A: Math. Gen. 34 (2001), 10377--10386.

\bibitem{I1} A. Iatrou. {\sl
Three dimensional integrable mappings,} arXiv:nlin.SI/0306052vl
  (2003)

\bibitem{I2} A. Iatrou. {\sl Higher dimensional integrable mappings,} Phys. D:
 179 (2003), 229--253.

\bibitem{IR1} A. Iatrou and , J.A.G. Roberts. {\sl Integrable mappings of the plane
preserving biquadratic invariant curves,} J. Phys. A: Math. Gen.
34 (2001), 6617--6636.

\bibitem{KYHRGO}K. Kimura, H. Yahagi, R. Hirota, A. Ramani, B. Grammaticos and Y. Ohta.
{\sl A new class of integrable discrete systems,} J. Phys. A:
Math. Gen. 35 (2002), 9205--9212.

\bibitem{K}  M.R.S. Kulenovi\'c. {\sl Invariants and related Liapunov
functions for difference equations}, Appl. Math. Lett. 13 (2000),
1--8.


\bibitem{L} S. Lang. ``Differential and Riemannian Manifolds''.
Springer. New York, 1995.

\bibitem{MT} J. Matsukidaira, D. Takahashi. {\sl Third--order
integrable difference equations generated by a pair of second
order equations}, J. Phys. A: Math. Gen. 39 (2006), 1151--1161.

\bibitem{RQ} J.A.G. Roberts, G.R.W. Quispel. {\sl Creating and
relating three dimensional integrable maps}, J. Phys. A: Math.
Gen. 39 (2006), L605--L615.

\bibitem{Z}  E. C. Zeeman. {\sl Geometric unfolding of a difference
equation}.
Unpublished paper. Hertford College, Oxford(1996), 1--42.
\end{thebibliography}
\end{document}